\documentclass[11pt,reqno]{amsart}
\usepackage{amsmath}
\usepackage{amssymb}
\usepackage{amsfonts}
\usepackage{graphicx}
\usepackage{amsthm}
\usepackage{enumerate}
\usepackage{lscape}
\usepackage{dsfont}
\usepackage{color}

\newcommand{\CH}{\mathds{C}\mathrm{H}}

\newcommand{\C}{\mathds{C}}

\newtheorem{theor}{Theorem}

\begin{document}
\title[On the epsilon function of Cartan--Hartogs domains]{A note on the coefficients of Rawnsley's epsilon function of Cartan--Hartogs domains}
\author[M. Zedda]{Michela Zedda}
\address{Dipartimento di Matematica ``G. Peano", Universit\`{a} di Torino}
\email{michela.zedda@gmail.com}

\thanks{
The author was supported by the project FIRB ``Geometria Differenziale e teoria geometrica delle funzioni'' and by INdAM-GNSAGA - Gruppo Nazionale per le Strutture Algebriche, Geometriche e le loro Applicazioni.
}
\date{}
 %32Q15  	Kähler manifolds. 53C55 Hermitian and Kählerian manifolds. 32A25: Integral representations; canonical kernels (Szeg?, Bergman, etc.) 32A07  	Special domains (Reinhardt, Hartogs, circular, tube) 32M15  	Hermitian symmetric spaces, bounded symmetric domains, Jordan algebras

\begin{abstract}
We  extend the result of Z. Feng and Z. Tu in \cite{fengtu} by showing that if one of the coefficients $a_j$, $2\leq j\leq n$, of  Rawnlsey's epsilon function associated to a $n$-dimensional Cartan--Hartogs domain is constant, then the domain is biholomorphically equivalent to the complex hyperbolic space. 
\end{abstract}

\maketitle

\section{Introduction and statement of main result}
Consider an $n$-dimensional complex manifold $(M,g)$ endowed with a K\"ahler metric $g$ and assume that there exists a globally defined  K\"ahler potential $\varphi\!:M\rightarrow \mathds{R}$ for $g$, i.e. if $\omega$ is the K\"ahler form associated to $g$, we have $\omega=\frac i2\partial\bar\partial \varphi$.
Let $\mathcal{H}_\alpha$ be the weighted Bergman space of square integrable holomorphic functions on $M$ with respect to the measure $e^{-\alpha\varphi}\frac{\omega^n}{n!}$, i.e.:
$$
\mathcal{H}_\alpha=\left\{f\in {\rm Hol(M)}|\,\int_Me^{-\alpha\varphi}|f|^2\frac{\omega^n}{n!}<\infty \right\}.
$$
If $\mathcal{H}_\alpha\neq \{0\}$, choose an orthonormal basis $\{f_j\}$ with respect to the product:
$$
(f,h)_\alpha=\int_Me^{-\alpha\varphi}f\bar h \frac{\omega^n}{n!},
$$
and denote by $K_\alpha(x,y)$ the reproducing kernel of $\mathcal{H}_\alpha$, namely:
$$
K_\alpha(x,y)=\sum_j f_j(x)\bar f_j(y),\qquad x,y\in M.
$$
Define the $\epsilon$-function associated to $g$ to be the function:
$$\epsilon_{\alpha g}(x)=e^{-\alpha\varphi(x)}K_{\alpha}(x,x), \quad x\in M.$$

In the literature the function $\epsilon_{\alpha g}$ was first introduced by J. Rawnsley under the name of $\eta$-{\em function} in \cite{rawnsley} and later as $\theta$-{\em function} in \cite{CGR}.

We say that  $\epsilon_{\alpha g}$ admits the {\em Engli\v{s} expansion}: 
\begin{equation}\label{englisexp}
\epsilon_{\alpha g}(x)\sim\sum_{j=0}^\infty a_j(x)\alpha^{n-j},\quad x\in M,
\end{equation}
for $\alpha\rightarrow+\infty$, if for every integers $l$, $r$ and every compact $H \subseteq M$,
\begin{equation}\label{TalphaMEANS}
\| \epsilon_{\alpha}(x) -
\sum_{j=0}^l  a_j(x){\alpha}^{n-j} \|_{C^r} \leq \frac{C(l, r, H)}{\alpha^{l+1}},
\end{equation}
for some constant $C(l, r, H) >0$. Such expansion is the counterpart for noncompact manifolds of the the celebrated TYZ (Tian-Yau-Zelditch) expansion
of Kempf's distortion function for polarized compact K\"ahler manifolds (see \cite{ze}
and  also \cite{arlquant}).
In \cite{englisasymp} M. Engli\v{s}
proved  that each of the coefficients $a_j(x)$ in \eqref{englisexp} 
 is a polynomial of the curvature of the metric $g$ and its covariant
derivatives at $x$, which can be
found by finitely many steps of algebraic operations, and gives an explicit expression of the coefficients $a_j$ for $j\leq 3$. \\

In this paper we consider the case of Cartan-Hartogs domains, which are defined as follows.
Consider a Cartan domain $\Omega\subset \C^d$, i.e. an irreducible bounded symmetric domain, of rank $r$ and numerical invariants $a$, $b$. Recall that the triple $\{r,a,b\}$ uniquely determines $\Omega$ and in particular it defines the dimension $d=\frac{r(r-1)}2a+rb+r$ and the genus $\gamma=(r-1)a+b+2$ of $\Omega$. Let $K(z, z)$ be the Bergman kernel of $\Omega$ and $N_\Omega(z, z)$ its {\em generic norm}, i.e.
\begin{equation}\label{genericnorm}
N_\Omega(z, z)=(V(\Omega)K(z, z))^{-\frac{1}{\gamma}},
\end{equation}
where $V(\Omega)$ is the total volume of $\Omega$ with respect to the Euclidean measure of the ambient complex Euclidean space.

For all positive real numbers $\mu$, a Cartan--Hartogs domains is given by $\left(M^{d_0}_{\Omega}(\mu),g(\mu)\right)$ where:
\begin{equation}\label{defm}
M^{d_0}_{\Omega}(\mu)=\left\{(z,w)\in \Omega\times\C^{d_0},\ ||w||^2<N_\Omega(z, z)^\mu\right\},
\end{equation}
and $g(\mu)$ is the K\"ahler metric whose associated K\"ahler form $\omega(\mu)$ can be described by the (globally defined) K\"ahler potential centered at the origin:
\begin{equation}\label{diastM}
\Phi(z,w)=-\log(N_{\Omega}(z, z)^\mu-||w||^2).
\end{equation}
The domain $\Omega$ is called  the {\em  base} of the Cartan--Hartogs domain 
$M^{d_0}_{\Omega}(\mu)$ (one also  says that 
$M^{d_0}_{\Omega}(\mu)$  is based on $\Omega$).
These domains have been considered by several authors (see e.g. \cite{roos} and references therein). In \cite{roos} it is shown that 
for $\mu_0=\gamma/(d+1)$, $(M^1_{\Omega}(\mu_0),g(\mu_0))$ is a complete K\"ahler-Einstein manifold which is homogeneous if and only if $\Omega$ is the complex hyperbolic space.
In \cite{articwall} the authors of the present paper proved that for $\Omega\neq\CH^d$, the metric $\alpha g(\mu)$ on $M^1_\Omega(\mu)$ is projectively induced for all positive real number $\alpha\geq \frac{(r-1)a}{2\mu}$, exhibing the first example of complete, noncompact, nonhomogeneous and projectively induced K\"ahler-Einstein metric. 
In \cite{ijgmmp} the  author of the present paper proved that for $d_0=1$, $g(\mu)$ is extremal (in the sense of Calabi \cite{Calabi82}) if and only if it is K\"ahler--Einstein and that the coefficient $a_2$ of Engli\v{s} expansion of the $\varepsilon$-function associated to a Cartan--Hartogs domain $(M^1_\Omega(\mu),g(\mu))$ is constant, then $(M^1_\Omega(\mu), g(\mu))$ is K\"ahler--Einstein, conjecturing also that the converse was true. In \cite{fengtu}, Z. Feng and Z. Tu generalize that theorem to generic $d_0$ and proved that conjecture. More precisely, they prove the following:

\begin{theor}[Z. Feng, Z. Tu {\cite[Th. 1.3]{fengtu}}]\label{fengzu}
The coefficient $a_2$ of the Rawnsley's $\epsilon$-function expansion is a constant on $M^{d_0}_\Omega(\mu)$ if and only if $\left(M^{d_0}_{\Omega}(\mu),g(\mu)\right)$ is biholomorphically isometric to the complex hyperbolic space $(\mathds{C}H^{d+d_0},g_{\rm hyp})$.
\end{theor}
Notice that  $g_{\rm hyp}$ denotes the hyperbolic metric on $\mathds{C}H^{d+d_0}$ and $$\left(\mathds{C}H^{d+d_0},g_{\rm hyp}\right)=\left(M^{d_0}_{\C H^{d}}(1),g(1)\right).$$

The prove of the previous theorem  is based on the  explicit formula for the $\epsilon$-function $\epsilon_{\alpha g(\mu)}$ of Cartan--Hartogs domains:
\begin{equation}\label{epsilonCH}
\epsilon_{\alpha g(\mu)}(z,w)=\frac 1{\mu^d}\sum_{k=0}^d\frac{D^k\tilde \chi(d)}{k!}\left(1-\frac{||w||^2}{N_\Omega(z, z)^\mu}\right)^{d-k}\frac{\Gamma(\alpha-d+k)}{\Gamma(\alpha-d-d_0)},
\end{equation}
for
$$
D^k\tilde \chi(d)=\sum_{j=0}^k{ k\choose j}(-1)^j\tilde \chi(d-j)
$$
and
$$
\tilde \chi(d-j)= \prod_{j=1}^r\frac{\Gamma(\mu (d-j)-\gamma-(j+1)\frac a2+2+b+ra)}{\Gamma(\mu (d-j)-\gamma+1+(j-1)\frac a2)},
$$
where $\Gamma$ is the usual $\Gamma$-function. Observe that formula \eqref{epsilonCH} shows that Engli\v{s} expansion of the $\epsilon$-function of Cartan--Hartogs domains is finite. In \cite{berezinCH} the author of this paper uses this formula to prove the existence of a Berezin-Engli\v{s} quantization for Cartan--Hartogs domains.

The aim of this paper is to generalize Theorem \ref{fengzu}  above (see next section) by proving the following: 
\begin{theor}\label{main}
For all $j=2,\dots, d+d_0$, any coefficient $a_j$ of the Rawnsley's $\epsilon$-function expansion is a constant on $M^{d_0}_\Omega(\mu)$ if and only if $\left(M^{d_0}_{\Omega}(\mu),g(\mu)\right)$ is biholomorphically isometric to the complex hyperbolic space $(\mathds{C}H^{d+d_0},g_{\rm hyp})$.
\end{theor}

The author would like to thank prof. Andrea Loi for the interesting and encouraging conversations.

\section{Proof of Theorem \ref{main}}

Due to Theorem \ref{fengzu} we need only to prove that if $a_j$ is constant for some $j=3,4,\dots, d+d_0$, then $a_2$ is.

Consider first the polynomial $P(\alpha)$ in the variable $\alpha$:
$$
P(\alpha)=\frac{\Gamma(\alpha-d+k)}{\Gamma(\alpha-d-d_0)},
$$
and observe that for $k=d$, ($d\geq 1$):
$$
P(\alpha)=\prod_{j=1}^{d+d_0}(\alpha-j),\quad \deg(P(\alpha))=d+d_0,%coefficienti dall'$a_0$ in poi});
$$
for $k=d-1$, ($d\geq 1$):
$$
P(\alpha)= \prod_{j=2}^{d+d_0}(\alpha-j),\quad  \deg(P(\alpha))=d+d_0-1,%coefficienti dall'$a_1$ in poi});
$$
for $k=d-2$, ($d\geq 2$):
$$
P(\alpha)=\prod_{j=3}^{d+d_0}(\alpha-j),\quad \deg(P(\alpha))=d+d_0-2,%, coefficienti dall'$a_2$ in poi});
$$
and so on.
Thus, from \eqref{epsilonCH}
we get that the factor with $k=d$ contributes to all the coefficients $a_0, a_1,\dots, a_{d+d_0}$ ($d\geq 1$), the factor with $k=d-1$ to all from $a_1$ to $a_{d+d_0}$ ($d\geq 1$), the factor with $k=d-2$ to all from $a_2$ ($d\geq 2$), and so on. Obviously the $j$-th coefficient is constant iff each one of its factors (except the $k=d$ one) vanishes, in fact the term $\left(1-\frac{||w||^2}{N_\Omega(z,\bar z)^\mu}\right)$ in each factor has a different power.

In particular, the coefficient $a_i$, $i=1,\dots, d+d_0,$ contains the factor:
$$
 \frac 1{\mu^d}\frac{D^{d-1}\tilde \chi(d)}{(d-1)!}\left(1-\frac{||w||^2}{N_\Omega(z,  z)^\mu}\right)A_{i}^2,\qquad (d\geq 1),
$$
and 
 the coefficient $a_i$, $i=2,\dots, d+d_0$, contains the factor:
$$
\frac 1{\mu^d}\frac{D^{d-2}\tilde \chi(d)}{(d-2)!}\left(1-\frac{||w||^2}{N_\Omega(z,  z)^\mu}\right)^{2}A_{i-1}^3,\qquad (d\geq 2),
$$
where we denote by $A_p^q$ the $p$-th coefficient of the polynomial in $\alpha$: 
$$
\prod_{j=q}^{d+d_0}(\alpha-j).
$$
Observe that $A_i^2$ and $A_{i-1}^3$ do not vanish. In fact we have:
$$\prod_{j=2}^{d+d_0}=\alpha^{d+d_0}+e_1(2,\dots,d+d_0)\alpha^{d+d_0-1}+\cdots+e_{d+d_0}(2,\dots,d+d_0),$$
$$\prod_{j=3}^{d+d_0}=\alpha^{d+d_0-1}+e_1(3,\dots,d+d_0)\alpha^{d+d_0-2}+\cdots+e_{d+d_0-1}(3,\dots,d+d_0)$$
where $e_j(x_1,\dots,x_n)$ is the elementary symmetric polynomial in the variables $(x_1,\dots,x_n)$, i.e.:
$$
e_j(x_1,\dots, x_n)=\sum_{1\leq k_1<k_2<\dots<k_j\leq n}x_{k_1}\dots x_{k_j}.
$$ 
Since in our case $x_j$ are positive integers, $A_i^2=e_i(2,\dots, d+d_0)$ and $A_{i-1}^3=e_{i-1}(3,\dots, d+d_0)$ do not vanish.

Thus we have that for $d\geq 2$ and for each $i=3,\dots, d$, if $a_i$ is constant then $D^{d-2}\tilde \chi(d)=D^{d-1}\tilde \chi(d)=0,$ and conclusion follows by \cite{fengtu}, where in the proof of Theorem 1.3 it is pointed out that when $d>1$ we have
$$D^{d-2}\tilde \chi(d)=D^{d-1}\tilde \chi(d)=0,$$
if and only if $a_2$ is constant.

If $d=1$, then $r=1$ and $\Omega=\mathds{C}H^1$, thus we need only to prove that if $a_j$ is constant for some $j=3,4,\dots, d_0+1$, then $\mu=1$. By the discussion above, if $a_j$ is constant for some $j=3,4,\dots, d_0+1$ then $D^{0}\tilde \chi(1)=0$, which by \cite[Lemma 3.5]{fengtu} directly implies $\mu=1$, concluding the proof.

\end{document}